\documentclass{article}

\usepackage{amsmath}
\usepackage{amsfonts}
\usepackage{amssymb}

\newtheorem{theorem}{Theorem}[section]
\newtheorem{prop}[theorem]{Postulate}

\date{}

\title{Calculation of Univariate Pad\'e Approximants \\
for solutions of the Michaelis-Menten equation with first order input using the $\tau$-method}
\author{Gareth Hegarty}
\begin{document}
\maketitle
\begin{abstract} 
In this paper the Jacobi formula is used to recursively generate (diagonal) univariate
Pad\'e approximants using the $\tau$-method 
for solutions Michaelis-Menten equation with first order input.
In the algorithm the Jacobi coefficients and error terms in the $\tau$-method 
are postulated to have a particular form, and this form is maintained by 
specific patterns of cancellations.
\end{abstract}

\section{Introduction}
The {\em Michaelis-Menten equation} originates in the field of 
pharmacokinetics where it models the elimination of a solute in a secondary compartment
(such as blood) over time. This equation may be written in the non-dimensional form
$$
\frac{dy}{dt}=\alpha e^{-t}-\frac{\beta y}{1+y},\ y(0)=0
$$
where $\alpha e^{-t}$ represents the first order input, $\alpha$ and $\beta$ are positive constants, 
and $y(t)$ is concentration at time $t$ (all non-dimensional). To do the asymptotic analysis 
we change the independent
variable to $x=e^{-t}$, converting the system to a singular ordinary differential of the form: 
\begin{equation}
x(1+y)y'+(\alpha x-\beta)y+\alpha x=0,\ y(0)=0\ (y(1)=0)
\label{mm1a}
\end{equation}
This equation is a special case of the second order {\em Riccati equation} 
(also known as the {\em Duffing equation} in physics) considered in Fair, Luke, \cite{fairluke}
which uses a series based approach to evaluate Pad\'e approximants using determinants. 
Here we use the $\tau$-method to evaluate the diagonal Pad\'e approximants
as this gives a more
detailed account of the approximants and errors in a manner similar to that used for the
first order Riccati equation in Fair, \cite{fair}. 

The general solution to (\ref{mm1a}) can be written as a bivariate series in $x$ and $x^\beta$
of the form:
$$
y=\sum^\infty_{n,m=0}c_{n,m}x^{n+m\beta},\ \ c_{0,0}=0
$$
However, the series does not generally converge in the interval $0\leq x\leq 1$, which is
needed to apply the conditon $y(1)=0$ (and hence solve for the unknown $c_{0,1}$).
To improve the approximation and enlarge the domain of convergence, one may also construct
bivariate Pad\'e approximants based on the bivariate series, as in Hegarty, \cite{hegarty}.
Though this is the broader aim of this work, here we focus on the univariate series in 
$x$ and $x^\beta$ separately as this is the foundation for the full bivariate case.
The series $y=y_{\infty,0}=\sum^\infty_{n=0}c_{n,0}x^n$ is required to satisfy the 
differential equation:
\begin{equation}
x(1+y)\frac{dy}{dx}+(\alpha x-\beta)y+\alpha x=0,\ y(0)=0\ y'(0)=
\alpha/(\beta-1)
\label{mm1b}
\end{equation}
Similarly, on the $z=x^\beta$ axis, the series $y=y_{0,\infty}=\sum^\infty_{m=0}c_{0,m}z^m$ 
is required to satisfy the differential equation:
\begin{equation}
z(1+y)\frac{dy}{dz}-y=0,\ y(0)=0\ y'(0)=c_1
\label{mm0}
\end{equation}
where $c_1=c_{0,1}$.
This may be solved using the {\em Lambert}-W function, which may also be written as a series:
$$y(x)=W_0(c_{0,1}x)=\sum^\infty_{k=1}\frac{(-k)^{k-1}(c_{0,1}x)^k}{k!}$$
which converges for $|x|<1/(c_{0,1}e)$. In this context we are more interested in 
how the Pad\'e approximants corresponding to this series may be evaluated.

\section{Formulation}
\label{s2}
Here we follow the $\tau$-method formulation in Fair, \cite{fair} to define 
the $n^{th}$ order main
diagonal Pad\'e approximant for the solution of (\ref{mm1b}) as a ratio of the form:
$$y_n=\frac{A_n}{B_n}=\frac{a^{(n)}_0+a^{(n)}_1x+...+a^{(n)}_nx^n}{1+b^{(n)}_1x+...+b^{(n)}_nx^n}$$
which satisfies the approximate equation:
\begin{equation}
B_n^3[x(1+y_n)y_n'+(\alpha x-\beta)y_n+\alpha x]=T_n=x^{2n}\sum^{n+1}_{k=1}\tau_{n,k}x^k
\label{Tndefn}
\end{equation}
This may also be written as:
$$
T_n=x(A_n+B_n)F_n+B_nH_n
$$
where
$$F_n=B_nA'_n-A_nB'_n,\ H_n=(\alpha x-\beta)A_nB_n+\alpha x B^2_n$$
To simplify the equations we also define the auxillary error term:
\begin{eqnarray}
U_n&=&x\{(A_n+B_n)L_n+(A_{n-1}+B_{n-1})F_n\}+B_nJ_n+B_{n-1}H_n\nonumber\\
&=&x^{2(n-1)}\sum^{n+2}_{k=1}u_{n,k}x^k
\label{Undefn}
\end{eqnarray}
where
$$L_n=A'_{n-1}B_n-B'_{n-1}A_n+A'_nB_{n-1}-B'_nA_{n-1},$$
$$J_n=(\alpha x-\beta)(A_{n-1}B_n+B_{n-1}A_n)+2\alpha x B_{n-1}B_n$$

The first couple of Pad\'e approximants can be calculated explicitly, e.g. 
starting with $n=0$, $A_0=0,\ B_0=1$, while for $n=1$ we get $A_1=-\alpha_1x$ 
($\alpha_1=-c_{1,0}=-\alpha/(-1+\beta)$, 
$B_1=1+\beta_1x$ ($\beta_1=-c_{2,0}/c_{1,0}=-\alpha\beta/((-1+\beta)(-2+\beta))$), 
where $c_{1,0}$, $c_{2,0}$ are the first couple of Taylor series coefficients around $x=0$.
Subsequent iterates can be generated using the Jacobi formulas:
\begin{eqnarray}
A_{n}=(1+\beta_n x) A_{n-1}+\alpha_n x^2 A_{n-2},\label{jaca}\\
B_{n}=(1+\beta_n x) B_{n-1}+\alpha_n x^2 B_{n-2},\label{jacb}
\end{eqnarray}
for $n>1$.
Recursion relations for $T_n$, $U_n$ can be derived from (\ref{jaca}), (\ref{jacb}), i.e.
expanding (\ref{Undefn}) gives:
\begin{eqnarray}
U_n &=& 3(1+\beta_n x)^2T_{n-1} + 3(1+\beta_{n-1} x)\alpha_n^2 x^4T_{n-2}\nonumber\\
&& + 2(1+\beta_n x)\alpha_n x^2 U_{n-1}
      + \alpha_n^2\alpha_{n-1}x^6 U_{n-2}\nonumber\\ 
&&+ 2\alpha_n^2\alpha_{n-1}x^6(A_{n-2}+B_{n-2})\Delta_{n-2}+ 2\alpha_n x^2(A_n+B_n)\Delta_{n-1}\nonumber\\
      &&+ (2+\beta_n x)\alpha_n x^2(A_{n-1}+B_{n-1})\Delta_{n-1}\label{Un0}
\end{eqnarray}
and expand (\ref{Tndefn}) as:
\begin{eqnarray}
T_n &=& (1+\beta_n x)^3T_{n-1} + \alpha_n^2x^4(\alpha_n x^2 + 3(1+\beta_n x)(1+\beta_{n-1} x))T_{n-2}\nonumber\\
&&      + (1+\beta_n x)^2\alpha_n x^2 U_{n-1}
      + (1+\beta_n x)\alpha_n^2\alpha_{n-1}x^6U_{n-2} \nonumber\\
      &&+ 2(1+\beta_n x)\alpha_n^2\alpha_{n-1}x^6(A_{n-2}+B_{n-2})\Delta_{n-2}\nonumber\\
      &&+ (2+\beta_n x)\alpha_n x^2(A_n + B_n)\Delta_{n-1}\label{Tn0}
\end{eqnarray}
where $\Delta_1=\alpha_1 x$ and $\Delta_n= -\alpha_nx^2\Delta_{n-1}=(-1)^{n-1}\alpha_{n,1}x^{2n-1};\ 
\alpha_{n,1}=\alpha_n...\alpha_1$ for $n>1$. 

From order $x^{2n-1}$ in $T_n$ we get the following formula for $\alpha_n$:
\begin{equation}
\alpha_n=\frac{-\tau_{n-1,1}}{u_{n-1,1}+2(-1)^{n-2}\alpha_{n-1,1}}
\label{alpha_n}
\end{equation}
and similarly, from order $x^{2n}$ we get $\beta_n$, i.e.
\begin{equation}
\beta_n=\frac{-\{\tau_{n-1,2} + \alpha_nu_{n-1,2} + 
2(a^{(n-1)}_1+b^{(n-1)}_1)(-1)^{n-2}\alpha_{n,1}\}}{
\{3\tau_{n-1,1} + 2\alpha_nu_{n-1,1} + 3(-1)^{n-2}\alpha_{n,1}\}}
\label{beta_n}
\end{equation}
Note that in the case $k=0$ we can use the formula for $\alpha_n$ above to get
\begin{equation}
u_{n,1}=\tau_{n-1,1}
\label{un1}
\end{equation}
To simplify the formula for $U_n$ we can also define the intermediary variable:
\begin{equation}
V_n=3(1+\beta_n x)T_{n-1}+\alpha_n x^2 U_{n-1}=x^{2(n-1)}\sum^{n+1}_{k=1}v_{n,k}x^k
\label{Vn}
\end{equation}
which allows us to rewrite the $U_n$ in (\ref{Un0}) as:
\begin{eqnarray}
U_n&=&3\beta_n^2x^2 T_{n-1}+2\beta_n\alpha_n x^3 U_{n-1} +\alpha_n^2 x^4V_{n-1}\nonumber\\
&&-3T_{n-1}+2V_{n}-(4+3\beta_n x)(A_{n-1}+B_{n-1})\Delta_n
\label{Una}
\end{eqnarray}
Similarly, using $V_n$ and $U_n$, we can rewrite the recursive formula for $T_n$ 
in (\ref{Tn0}) as:
\begin{eqnarray}
T_n&=&\beta_n^3x^3 T_{n-1}+\beta_n^2\alpha_n x^4 U_{n-1} +\beta_n\alpha_n^2 x^5V_{n-1}
+\alpha_n^3 x^6T_{n-2}
\nonumber\\
&&+U_n-V_{n}+T_{n-1}
+(2+\beta_n x)(A_{n-1}+B_{n-1})\Delta_n\nonumber\\
&&+2\beta_n\alpha_n x^3(A_{n-2}+B_{n-2})\Delta_n
-\beta_n x(A_{n}+B_{n})\Delta_n\label{Tna}
\end{eqnarray}

\section{Recursive Algorithm}
\label{s3}
While the formulas above may be used to calculate the Jacobi coefficents $\alpha_n$, $\beta_n$
and error coefficients $v_{n,k}$, $u_{n,k+1}$, $\tau_{n,k}$ directly,
there is an inherent structure in these parameters which can be used to simplify these
calculations. 
To characterise this structure we introduce monomials $P_n$, $Q_n$ and 
$M_{n,k}$, $J_{n,k}$, $L_{n,k+1}$ and $K_{n,k}$
which are all made up of factors $(-p+\beta)^*$
for $p\in{\cal C}^{(n)}=\{2,...,2n\}$. 
We can group the integers $p\in{\cal C}^{(n)}$
into 4 distinct subsets, 
i.e. ${\cal C}^{(n)}={\fam2 C}^{(n)}_0\cup{\fam2 C}^{(n)}_1\cup
{\fam2 C}^{(n)}_2\cup{\fam2 C}^{(n)}_3$ where:
$${\fam2 C}^{(n)}_1=\{p\in {\cal C}^{(n)}:\hbox{Mod}[n,p]=0\}$$
$${\fam2 C}^{(n)}_2=\{p\in {\cal C}^{(n)}:p<2(\hbox{Mod}[n,p]-1)\}$$
$${\fam2 C}^{(n)}_3=\{p\in {\cal C}^{(n)}:p=2(\hbox{Mod}[n,p]-1)\}$$
$${\fam2 C}^{(n)}_0=\{p\in {\cal C}^{(n)}:p>2(\hbox{Mod}[n,p]-1)\ \hbox{and}\ \hbox{Mod}[n,p]\neq 0\}$$ 
Further, one can partition the sets
${\cal C}^{(n)}_0$ and ${\cal C}^{(n)}_2$ into subsets ${\cal C}^{(n)}_0[r]$,
${\cal C}^{(n)}_2[r]$ respectively for $r>0$, i.e. $p\in{\cal C}^{(n)}_0[r]$ 
or $p\in{\cal C}^{(n)}_2[r]$ if
$$r=|p-2(\hbox{Mod}[n,p]-1)|$$
We also introduce coefficients $p_{n,1}$, $p_{n,2}$ shown in table \ref{tab1}
which may be used to define the mononials:
\begin{equation}
P_n=-p_{n,1}p_{n-1,2}^2p_{n-1,1},\ Q_n=p_{n,2}p_{n-1,2};\ n>1
\label{PnQn}
\end{equation}
along with $P_1=-1$, $Q_1=(-2+\beta)$.
These monomials will appear in the formulas for $\alpha_n$ and $\beta_n$ respectively 
in Postulate~\ref{prop1}.

The coefficients $p_{n,1}$, $p_{n,2}$ may also be used to define the first (explicitly)
and last (recursively) cases for $M_{n,k}$:
\begin{equation}
M_{n,1}=p_{n,2},\ M_{n,n}=p_{n,1}p_{n,2}M_{n-1,n-1}
\label{Mn1Mnn}
\end{equation}
In between these values we expect the power in $M_{n,k}:(-p+\beta)^*$ to be non-decreasing 
with $k$, and increasing by at most $1$ for each step.
To characterise these increases, we denote the $k$ at which $s^{th}$ 
increase in $M_{n,k}$ by $k=k^{M,s}_{n,p}:(-p+\beta)^+\ (\leq n)$.
In the case $k=1$ this may be given explicitly as:
\begin{equation}
k=k^{M,1}_{n,p}=1+2n_p-p+2\hbox{Mod}[n-n_p,p],\ n_p=\hbox{ceil}(p/2)
\label{kM1}
\end{equation}
While this formula is true
for any $p\in{\cal C}^{(n)}$, we can simplify this in particular cases, e.g.
\begin{equation}
k^{M,1}_{n,p}=\left\{\begin{array}{ll}
r+3,&p\in{\cal C}^{(n)}_2[r],\\
3,&p\in{\cal C}^{(n)}_3\\
2,&p\in{\cal C}^{(n)}_0[1]\\
1,&p\in{\cal C}^{(n)}_0[2]\\
p,&p\in{\cal C}^{(n)}_1\cup{\cal C}^{(n)}_0[>2]
\end{array}\right.
\label{kM1np}
\end{equation}
Beyond this point the power increases periodically (every $p$ steps in $k$) until 
it reaches the value in $M_{n,n}$, i.e.
$$
k^{M,s}_{n,p}=k^{M,1}_{n,p}+(s-1)p\leq n
$$
It is also obvious that $k^{M,s}_{n,p}$ is $p$-periodic in $n$ and $k$, i.e.
$$
k^{M,s}_{n+p,p}=k^{M,s}_{n,p},\ k^{M,s+1}_{n,p}=k^{M,s}_{n,p}+p;\ (\leq n)
$$

We introduce the monomials $J_{n,k}$, $L_{n,k}$, $K_{n,k}$ together since these are 
inter-related.
The first two $J_{n,k}$ are given explicitly by:
$$
J_{n,1}=(-1)^n(P_n...P_1),\ 
J_{n,2}=Q_n(-1)^n(P_n...P_1)
$$
while the last two may be defined recursively:
\begin{equation}
J_{n,n+1}=p_{n,2}p_{n,1}K_{n-1,n},\ 
J_{n,n}=p_{n,2}p_{n,1}(-n+\beta)K_{n-1,n-1}
\label{Jnnp1etc}
\end{equation}
Similarly, the first two $L_{n,k}$ are given explicitly by:
$$
L_{n,1}=\frac{(-1)^n(P_n...P_1)}{(-(2n-1)+\beta)},\ 
L_{n,2}=p_{n,2}(-1)^n(P_n...P_1),
$$
while the last two may be defined recursively:
\begin{equation}
L_{n,n+2}=p_{n,2}p_{n,1}J_{n,n+1},\ L_{n,n+1}=p_{n,2}p_{n,1}J_{n,n}
\label{LnnpLn2np1}
\end{equation}
Given the equation (\ref{un1}), we define\footnote{
Another simplification which could be made is to take the common factor
of $(-1)^n(P_n...P_1)$ out of $J_{n,k}$, $L_{n,k+1}$ and $K_{n-1,k}$ (noting that they all start 
with this factor for $k=1$), i.e.
$$J_{n,k}=(-1)^n(P_n...P_1)J'_{n,k},\ 
L_{n,k+1}=(-1)^n(P_n...P_1)L'_{n,k+1},\ 
K_{n,k}=(-1)^n(P_n...P_1)K'_{n,k}$$
etc.} 
the first two $K_{n-1,k}$:
$$
K_{n-1,1}=L_{n,1},\ 
K_{n-1,2}=\frac{p_{n-1,2}L_{n,2}}{(-2n+\beta)(-(2n-1)+\beta)}
$$
while the last two may be defined recursively:
\begin{equation}
K_{n,n+1}=p_{n,2}p_{n,1}L_{n,n+2},\ K_{n,n}=p_{n,2}p_{n,1}L_{n,n+1}/(-(n+1)+\beta)
\label{Knnp1Knn}
\end{equation}
Similar to what we did for $M_{n,k}$, we now consider the patterns of increases
in $J_{n,k}, L_{n,k}, K_{n,k}:(-p+\beta)^+$ as $k$ increases. 
To characterise these increases, we denote the $k$ at which $s^{th}$ 
increase in $J_{n,k}$, $L_{n,k}$, $K_{n,k}$ occur by
$$k^{J,s}_{n,p},\ k^{L,s}_{n,p},\ k^{K,s}_{n,p}:(-p+\beta)^+$$
In between and consistent with these end values we can determine the powers of $(-p+\beta)^*$
in the monomials by defining the $k$ at which the increases occur
for each case $p\in{\cal C}^{(n)}$ separately.
For $p\in {\cal C}^{(n)}_3$ we have:
$$k^{J,s}_{n,p}=\left\{\begin{array}{ll}
2,&s=1\\
5+(s-2)p,&s\geq2
\end{array}\right.$$
For $p\in {\cal C}^{(n)}_2[r]$ we have:
$$k^{J,s}_{n,p}=\left\{\begin{array}{ll}
r+2,&s=1\\
2r+5+(s-2)p,&s\geq 2
\end{array}\right.$$
We can them use this to define the $L_{n,k}$ and $K_{n,k}$ increases in 
both of the previous cases:
$$
k^{L,s}_{n,p}=k^{K,s}_{n,p}=k^{J,s}_{n,p}+2,\ p\in{\cal C}^{(n)}_3\cup{\cal C}^{(n)}_2
$$ 
For $p\in {\cal C}^{(n)}_1$ we have:
$$k^{J,s}_{n,p}=sp\ (\leq n),\ k^{K,s}_{n,p}=k^{J,s}_{n,p},\ k^{L,1}_{n,p}=k^{J,s}_{n,p}+1$$ 
Next consider $p\in {\cal C}^{(n)}_0[2]$ where we see similar increases for $J_{n,k}$
and $K_{n,k}$ at $k=2$ and $k=3$, i.e.
$$
k^{J,s}_{n,p}=k^{K,s}_{n,p}=\left\{\begin{array}{ll}
2,&s=1\\
3+(s-2)p,&s\geq2
\end{array}\right.
$$
while $L_{n,k}$ increases three times at $k=2$, $k=3$
and $k=4$, i.e.
$$
k^{L,1}_{n,p}=2,\ k^{L,s}_{n,p}=k^{J,s}_{n,p}+1;\ s>1
$$
Typical of this is the example $n=6,p=4$ where the increases are shown 
in table \ref{tabknp}. 
\begin{table}
\centering
\begin{tabular}{|c|c|c|c|}
\hline
$k$ & $J_{n,k}$ & $L_{n,k}$ & $K_{n,k}$\\
\hline
$4$ & $(-4+\beta)^4$ & $(-4+\beta)^5$ & $(-4+\beta)^6$\\
$3$ & $(-4+\beta)^4$ & $(-4+\beta)^4$ & $(-4+\beta)^6$\\
$2$ & $(-4+\beta)^3$ & $(-4+\beta)^3$ & $(-4+\beta)^5$\\
$1$ & $(-4+\beta)^2$ & $(-4+\beta)^2$ & $(-4+\beta)^4$\\
\hline
\end{tabular} 
\caption{Initial increases in $J_{n,k}, L_{n,k}, K_{n,k}$ for $n=6,p=4\in{\cal C}^{(6)}_0[2]$}
\label{tabknp}
\end{table}
For $p\in {\cal C}^{(n)}_0[1]$ we see increases in $L_{n,k}$ and $K_{n,k}$
at $k=3$ and $k=5$ and after this we get periodic increases until reaching
the respective terminal exponent, while $J_{n,k}$ increases only once at $k=3$ 
before increasing periodically, i.e.
$$k^{J,s}_{n,p}=3+(s-1)p\ (\leq n)$$
and
$$k^{L,1}_{n,p}=k^{K,1}_{n,p}=3,\ k^{L,s}_{n,p}=k^{K,s}_{n,p}=k^{J,s-1}_{n,p}+2;\ s>1$$

\begin{prop}
\label{prop1}
The Jacobi coefficients in (\ref{alpha_n}), (\ref{beta_n}) may be written explicitly 
in the form:
\begin{equation}
\alpha_n=\frac{c_1^2\chi_{n-2}\chi_n}{P_n\chi^2_{n-1}},\ \beta_n
=\frac{c_1\psi_n}{Q_n\chi_{n-1}\chi_n},\ n\geq 1
\label{alpnbetn}
\end{equation}
where $P_n$, $Q_n$ are given in (\ref{PnQn}) and $\chi_{-1}=\beta^2/\alpha$, 
$\chi_0=\beta$, $\chi_n=S_{n-1,1}$ for $n>0$, and
\begin{eqnarray*}
\psi_n&=&-\big\{S_{n-1,2}(-2n+\beta)(-(2n-1)+\beta)(-1+\beta)^{\delta_{n,2}}
\\
&&-p_{n,2}\chi_nW_{n-1,2}(-1+\beta)^{\delta_{n,2}}
+2p_{n,2}\chi_n\varphi^{(n-1)}_1\big\}/(-2n+\beta),\ n>1
\end{eqnarray*}
We can also write the components of $A_n+B_n$ in the form:
\begin{equation}
a^{(n)}_0+b^{(n)}_0=1,\ a^{(n)}_k+b^{(n)}_k=\frac{\alpha^k\varphi^{(n)}_k}{M_{n,k}\chi_n};\ \ 0<k\leq n
\label{phin}\nonumber
\end{equation}
where $\varphi^{(n)}_k$ are polynomials in $\beta$. The errors also have explicit formulas:
\begin{equation} 
v_{n,k}=\frac{c_1^{2(n-1)+k}\beta G_{n,k}}{J_{n,k}\chi_{n-1}^{\min\{k,2\}}\chi_n^{\min\{k-2,1\}}},\ 
k=1,...,n+1,\ n>0
\nonumber
\end{equation}
\begin{equation} 
u_{n,k}=\frac{c_1^{2(n-1)+k}\beta W_{n,k}}{L_{n,k}\chi_{n-1}\chi_n^{\min\{k-2,2\}}},\ 
k=1,...,n+2,\ n\geq 0
\nonumber
\end{equation}
and
\begin{equation}
\tau_{n,k}=\frac{c_1^{2n+k}\beta S_{n,k}}{K_{n,k}\chi^{\min\{k,3\}}_n},\ 
k=1,...,n+1,\ n\geq 0
\nonumber
\end{equation}
where $S_{n,k}$, $W_{n,k+1}$ and $G_{n,k}$ are polynomials in $\beta$
In addition, in the last few error terms the numerators can be factorised, i.e.
\begin{equation}
G_{n,n+1}=\beta\Omega_{n-1}G'_{n,n+1},\ W_{n,n+2}=\beta\Omega_nW'_{n,n+2}\nonumber
\end{equation}
\begin{equation}
S_{n,n}=\Omega_nS'_{n,n},\ S_{n,n+1}=\beta\Omega^2_nS'_{n,n+1},\nonumber
\end{equation}
where $\Omega_n$ are also polynomials in $\beta$.
\end{prop}

We can verify these expressions for the first few $n$ by using the expressions from section \ref{s2},
e.g. in the case $n=0$:
$$\Omega_0=1,\ W_{0,1}=1+\beta,\ L_{0,1}=(-1+\beta), W_{0,2}=-1,\ L_{0,2}=1,$$
and 
$$S_{0,1}=1,\ K_{0,1}=(-1+\beta)^{-1}=L_{1,1},
\varphi^{(0)}_0=\beta,\ 
M_{0,0}=1$$
Similarly, for $n=1$ we have $\Omega_1=1$ and:
$$G_{1,1}=2(-1+\beta),\ J_{1,1}=1,\ G_{1,2}=-2(1+\beta),\ J_{1,2}=(-2+\beta)(-1+\beta)^{-1},$$
and since $u_{1,1}=\tau_{0,1}=\alpha$, we get $W_{1,1}=1$ and similarly:
$$W_{1,2}=-2(1+\beta),\ L_{1,2}=(-2+\beta),\ W'_{1,3}=4+\beta,\ L_{1,3}=(-2+\beta)^2$$
and
$$S'_{1,1}=-4+3\beta,\ K_{1,1}=(-2+\beta)^2,\ S'_{1,2}=\varphi^{(1)}_1=-2,\ 
K_{1,2}=(-2+\beta)^3,$$ 
$$\varphi^{(1)}_0=0,\ M_{1,0}=1,\ \varphi^{(1)}_1=-2,\ M_{1,1}=(-2+\beta)$$
In the following subsections we calculate the numerators $\varphi^{(n)}_k$, $G_{n,k}$, 
$W_{n,k+1}$ and $S_{n,k}$ in the general case ($n>1$) and show how these remain
polynomials due to patterns of cancellations with respect to both the $\chi_{n-1}^*$
and $(-p+\beta)^*$ appearing in the denominators.
Since the monomials are made up of factors $(-p+\beta)^*$ for $p=2,...,2n$, we can focus on each 
factor $(-p+\beta)^*$ separately and distinguish between {\em major factors} ($p\geq n$)
and {\em minor factors} ($p<n$). 
In each case we find an obvious (diagonal) pattern of increases (and cancellations) 
in the powers of the major factors
which is then propagated to the minor factors using periodicity over $n$ and $k$. 

{\noindent\bf Remark.} Consider the example of the equation (\ref{mm1b}) in the 
case $\alpha=c^2_1$, $\beta=1+c_1$, and hence $y'(0)=\alpha/(\beta-1)=c_1$.
It is clear that for small $c_1\sim 0$ the solution of this equation will approximate
the solution of (\ref{mm0}). Moreover,  we can use this to evaluate the Jacobi coefficients
$\alpha_n$, $\beta_n$ for (\ref{mm0}) as:
$$\alpha_n=c_1^2\lim_{c_1\to 0}\frac{\chi_{n-2}\chi_n}{P_n\chi^2_{n-1}},\ \beta_n
=c_1\lim_{c_1\to 0}\frac{\psi_n}{Q_n\chi_{n-1}\chi_n}$$
and similarly for the error terms. One can also improve on this to account for the 
obvious cancellations as one applies the limit as $c_1\to 0$.
\subsection{Calculation of $a^{(n)}_k+b^{(n)}_k$}
\label{ssMalg}
We consider the calculation of $a^{(n)}_k+b^{(n)}_k$ first as is of a simpler form and can 
be solved independently from the others which follow.
We start by
considering the components $x^{k}:A_n+B_n$ and using the Jacobi identity and the Postulate 
to expand. In the case $k=1$ we can use (\ref{Mn1Mnn}) to get:
$$
p_{n-1,2}\chi_{n-1}\varphi^{(n)}_1=p_{n,2}\chi_n\varphi^{(n-1)}_1+\psi_n
$$
Similarly, in the case $k=n$ the recursive formula for $M_{n,n}$
in (\ref{Mn1Mnn}) allows us to derive the expression:
\begin{equation}
p_{n-1,2}\chi^2_{n-1}\varphi^{(n)}_n=p_{n,1}\psi_n\varphi^{(n-1)}_{n-1} - 
p_{n,2}\chi_n^2\varphi^{(n-2)}_{n-2}
\label{phinneqn}
\end{equation}
In the general case $1<k<n$ we have:
\begin{equation}
\frac{\chi_{n-1}^2\varphi^{(n)}_{k}}{M_{n,k}}
=\frac{\chi_{n-1}\chi_n\varphi^{(n-1)}_{k}}{M_{n-1,k}}+
\frac{\varphi^{(n-1)}_{k-1}\psi_n}{Q_nM_{n-1,k-1}}+\frac{\varphi^{(n-2)}_{k-2}\chi_n^2}{P_nM_{n-2,k-2}}
\label{abk}
\end{equation}
Thus we can write the polynomial $\varphi^{(n)}_{k}$ as a quotient:
\begin{equation}
\varphi^{(n)}_{k}=\big\{c_1\chi_{n-1}\chi_n\varphi^{(n-1)}_{k}
+c_2\varphi^{(n-1)}_{k-1}\psi_n+c_3\varphi^{(n-2)}_{k-2}\chi_n^2\big\}/(c\chi_{n-1}^2)
\label{phinklcm}
\end{equation}
where:
$$
c_1=\frac{LCM_{n,k}}{M_{n-1,k}},\ c_2=\frac{LCM_{n,k}}{Q_nM_{n-1,k-1}},\ 
c_3=\frac{LCM_{n,k}}{P_nM_{n-2,k-2}},\ 
c=\frac{LCM_{n,k}}{M_{n,k}}
$$
and
\begin{equation}
\hbox{LCM}_{n,k}=
\hbox{LCM}\{M_{n-1,k},Q_{n}M_{n-1,k-1},P_{n}M_{n-2,k-2}\}
\label{lcmM}
\end{equation}
This leads us to focus on the cancelled factors in $c=\hbox{LCM}_{n,k}/M_{n,k}$.

The first thing one notices is the
obvious pattern in the cancellation of major factors $\hbox{LCM}_{n,k}\to M_{n,k}:(-p+\beta)^*$
for $n\leq p\leq 2(n-1)$.
These cancellations generally occur where the LCM increases and 
come in 4 {\em types} according to the relative orders of $(-p+\beta)^*$
in the three components of $\hbox{LCM}_{n,k}^M$, as in table \ref{tabctypes}.
\begin{table}
\centering
\begin{tabular}{|c|c|c|c|}
\hline
$(-p+\beta)^*$ & $M_{n-1,k}$ & $Q_nM_{n-1,k-1}$ & $P_nM_{n-2,k-2}$\\
\hline
$[1]$ & $+0$ & $+1$ & $+1$\\
$[2]$ & $+1$ & $+1$ & $+1$\\
$[3]$ & $+1$ & $+0$ & $+1$\\
$[4]$ & $+1$ & $+1$ & $+0$\\
\hline
\end{tabular} 
\caption{Cancel types in $M_{n,k}$}
\label{tabctypes}
\end{table}
Typical of this pattern is the $n=5$ case where we see:
$$M_{5,5}=(-10+\beta)(-9+\beta)(-8+\beta)[1](-7+\beta)(-6+\beta)(-5+\beta)...$$
$$M_{5,4}=(-10+\beta)(-9+\beta)(-8+\beta)[1](-7+\beta)(-6+\beta)^0[2](-5+\beta)^0[3]...$$
$$M_{5,3}=(-10+\beta)(-9+\beta)(-8+\beta)[1](-7+\beta)^0[2](-6+\beta)^0[3](-5+\beta)^0...$$
$$M_{5,2}=(-10+\beta)(-9+\beta)(-8+\beta)^0[1,2](-7+\beta)^0[3](-6+\beta)^0(-5+\beta)^0...$$
$$M_{5,1}=(-10+\beta)(-9+\beta)^0(-8+\beta)^0[4](-7+\beta)^0(-6+\beta)^0(-5+\beta)^0...$$
For $p=n\in{\cal C}^{(n)}_1$ we see a single cancellation $[3]$
with the increase in LCM at $k=n-1$, followed by the increase at $k=n$ (consistent
with (\ref{kM1np})). 
Similarly, for 
$n<p<2(n-1)$ (i.e. $p\in{\cal C}^{(n)}_2[r]$ where $r=2(Mod[n,p]-1)-p$) we see a diagonal pattern of 
consecutive cancellations $[3],[2]$ at each increase in the LCM.
Finally, for $p=2(n-1)\in{\cal C}^{(n)}_3$ the $M_{n,k}:(-p+\beta)^*$ cancels at every $k$ and twice
$[1,2]$ at each increase in the LCM (with the power increases after this).
The second point to notice is that these cancellations in the major factors are repeated periodically
(in the minor factors), every $p$ steps in $n$ and $k$.
For instance, in the $n=5$ example above, the
single cancellation $M_{5,4}:(-5+\beta)^0[3]$ will be projected to $M_{10,4},M_{10,9}:(-5+\beta)^0[3]$
(and similarly to $M_{15,4},M_{15,9},M_{15,14}:(-5+\beta)^0[3]$,.. etc.), 
while the pair $M_{5,3},M_{5,4}:(-6+\beta)^0[*]$ will be projected to the pair
$M_{11,3},M_{11,4}:(-6+\beta)^0[*]$ and again at $M_{11,9},M_{11,10}:(-6+\beta)^0[*]$ and so on.

This cancellation pattern can easily be extended to all of $p\in{\cal C}^{(n)}$,
including $p\in{\cal C}^{(n)}_0$ where $M_{n,k}:(-p+\beta)^*$ just increases with the LCM
without cancellation. Moreover, since these cancellation must be consistent with the increases
seen in (\ref{kM1np}), $\hbox{LCM}_{n,k}:(-p+\beta)^*$ must increase at $k^{M,s}_{n,p}-1$
for $p\in{\cal C}^{(n)}_1\cup {\cal C}^{(n)}_3$, at $k^{M,s}_{n,p}-2$
for $p\in{\cal C}^{(n)}_2$, and at $k^{M,s}_{n,p}$ for $p\in{\cal C}^{(n)}_0$.
In all cases, it is easy to justify this conditions are met
pattern in general before considering edge cases. For example, in the case $p\in{\cal C}^{(n)}_1$, 
in the general case $n>6$ we get 
$p\in{\cal C}^{(n-1)}_2[n-4]\cap{\cal C}^{(n-2)}_2[n-6]$ and hence for $k=n-1$
we get $M_{n-1,n-1},M_{n-2,n-3}:(-p+\beta)^+$, justifying the $[3]$ cancellation pattern.

\subsection{Calculation of $v_{n,k}$}
\label{ssJalg}
To simplify the calculation of $v_{n,k}$ (and deal with edge cases), it is helpful 
work in terms of {\em modified} numerators:
$$
v_{n,k}=\frac{c_1^{2(n-1)+k}\beta \tilde G_{n,k}}{J_{n,k}\chi_{n-1}^{2}\chi_n^{1}},\ 
\tilde G_{n,k}=\chi_{n-1}^{2-\min\{k,2\}}\chi_n^{1-\min\{k-2,1\}}G_{n,k},
$$
Expanding $x^k:V_n$ from (\ref{Vn}) we get a similar equation to (\ref{abk}), i.e.
\begin{equation}
\frac{\chi_{n-1}^2\tilde G^{(n)}_{k}}{J_{n,k}}=
\bigg(\frac{3\tilde S_{n-1,k}
\chi_{n-1}\chi_{n}}{K_{n-1,k}}\bigg)+\bigg[\frac{3\tilde S_{n-1,k-1}\psi_{n}
}{Q_{n}K_{n-1,k-1}}
+\frac{\tilde W_{n-1,k}
\chi_{n}^2}{P_{n}L_{n-1,k}}\bigg]
\label{gjnm1}
\end{equation}
In the case $k=1$ this simplifies to:
$$
G_{n,1}=2(-2n+\beta),
$$
The equation in (\ref{gjnm1}) may similarly be expressed in terms of LCMs as in (\ref{phinklcm}), 
leading us to focus on the ratios:
$$
\hbox{LCM}\{K_{n-1,k},Q_{n}K_{n-1,k-1},P_{n}L_{n-1,k}\}/J_{n,k}
$$
as these represent the cancelled factors. Given this similarity, it is not surprising
that we get a similar pattern of cancellations to those we found in $M_{n,k}$ above (with 
some minor differences). The most obvious difference is that we get no cancellations of
in the major factors, apart from for $p=2(n-1)\in {\cal C}^{(n)}_3$.
The factors $p\in {\fam2 C}^{(n)}_3$ have the same behaviour as in $M_{n,k}$,
cancelling twice at an increase in the LCM or once otherwise (for $k>2$).
In between these extremes (i.e. $3\leq k\leq n-1$) one does see a similar pattern of 
increases in LCMs for adjacent pairs of {\em major factors} ($p\geq n$)
but in this case these do not cancel. These increases in the major factors are also projected
to increases in the minor factors, i.e. $J_{n,k}:(-p+\beta)^+\to J_{n+p,k}:(-p+\beta)^+$, 
essentially meaning that minor factors in ${\fam2 C}^{(n)}_1$ or ${\fam2 C}^{(n)}_2$ 
do not cancel 
at the first increase in the $\hbox{LCM}_{n,k}$. However, beyond this point 
(i.e. for $n<k\leq n+p$) they cancel as before, i.e. once for $(-p+\beta)\in {\fam2 C}^{(n)}_1$
and as a consecutive pair for $(-p+\beta)\in {\fam2 C}^{(n)}_2$.

In the case $k=n+1$ 
we can use (\ref{Jnnp1etc}) and (\ref{Knnp1Knn}) to simplify (\ref{gjnm1}):
\begin{equation}
p_{n-1,2}\chi_{n-1}^2\tilde G_{n,n+1}=
3p_{n,1}\psi_n\tilde S_{n-1,n}-p_{n,2}\chi_n^2\tilde W_{n-1,n+1}
\label{Gnnp1imp}
\end{equation}
While this can be used to factorise $\tilde G_{n,n+1}=\beta\Omega_{n-1}\tilde G'_{n,n+1}$,
we can also get an explicit formula for $\tilde G'_{n,n+1}$ itself:
\begin{equation}
\tilde G'_{n,n+1}=-(3\Omega_{n-1}\tilde S'_{n,n+1}+2\tilde W'_{n,n+2})
\label{Gdnnp1}
\end{equation}

\subsection{Calculation of $u_{n,k+1}$}
\label{ssLalg}
To simplify the calculation of $u_{n,k+1}$ (and deal with edge cases), it is helpful 
work in terms of {\em modified} numerators:
$$
u_{n,k+1}=\frac{c_1^{2(n-1)+k+1}\beta \tilde W_{n,k+1}}{L_{n,k+1}\chi_{n-1}\chi_n^{2}},\ 
\tilde W_{n,k+1}=\chi_{n}^{2-\min\{k-1,2\}}W_{n,k+1},
$$
In order to cancel the $\chi_{n-2}^2$, we note the implicit formula:
\begin{eqnarray}
3p_{n,1}\Omega_{n-1}\psi_{n}S'_{n-1,n}-p_{n,2}\chi^2_{n}W'_{n-1,n+1}\ \ \ \ \ \ \ \ \ \ \ \ \ &&\nonumber\\
= -p_{n-1,2}\chi^2_{n-1}(3\Omega_{n-1}S'_{n,n+1}+2W'_{n,n+2})&&
\label{aaa}
\end{eqnarray}
Thus, expanding $x^{2(n-1)+k+1}:U_n$ from (\ref{Una}) and using the postulate we find:
$$
u_{n,k+1}=\frac{\alpha^{2n+k-1}\beta}{\chi_{n-1}\chi^2_{n}}.
\frac{1}{\chi^4_{n-1}}\bigg\{
\bigg(\frac{3\psi_n^2\tilde S_{n-1,k-1}
}{Q_n^2K_{n-1,k-1}}\bigg)
+\bigg(\frac{2\psi_n\tilde W_{n-1,k}
\chi_n^{2}}{Q_nP_nL_{n-1,k}}\bigg)
$$
$$
+ \bigg(\frac{\chi_{n}^4\tilde G^{(n-1)}_{k-1}}{P_n^2J_{n-1,k-1}}\bigg)
-\bigg(\frac{3\tilde S_{n-1,k+1}\chi_{n-1}^{2}\chi_n^{2}}{K_{n-1,k+1}}\bigg)
+ \bigg(\frac{2\chi_{n}^4\chi_{n-1}^3\tilde G^{(n)}_{k+1}}{J_{n,k+1}}\bigg)
$$
$$
+3(-1)^{n-2}\bigg(\frac{\psi_n\varphi^{(n-1)}_{k-1}
\chi_{n-1}^2\chi_n^{2}}{Q_n(P_n...P_1)M_{n-1,k-1}}\bigg)
+4(-1)^{n-2}\bigg(\frac{\varphi^{(n-1)}_{k}
\chi_{n-1}^3\chi_n^{3}}{(P_n...P_1)M_{n-1,k}}\bigg)
\bigg\}
$$
Labelling the components inside the brackets as (I)-(VII) for reference, 
we can break this sum into two steps. The first step\footnote{Note that this step can actually 
be decomposed further, i.e. adding (I)+(II)+(III) to get a $\chi_{n-1}^2$ term before adding 
(IV) and (VI) to get to $\chi_{n-1}^3$. Though this would be useful improving the algorithm efficiency,
it is not essential to the argument here.} is to add all terms apart from (V)
and (VII) to get a $\chi_{n-1}^3$ term: 
$$
(I) + (II) + (III) +(IV) + (VI)
=\frac{\chi_{n-1}^3X^{(3)}_{n,k+1}}{R^{(3)}_{n,k+1}}$$
where $R^{(3)}_{n,k+1}$ is determined from the LCM of denominators in this sum,
reducing the power of $(-p+\beta)^*$ in the LCM to match the leading term in 
the denominators of (V) and (VII) for $p\in {\cal C}^{(n)}_1$, ${\cal C}^{(n)}_2$
and ${\cal C}^{(n)}_3$, i.e. those in $J_{n,k+1}$ according to (\ref{JgeqM}).
An important consequence of these patterns of increases that is used in 
the algorithm for $L_{n,k+1}$ is that for any 
$p\in {\cal C}^{(n)}_1\cup {\cal C}^{(n)}_2\cup {\cal C}^{(n)}_3$ and
any $k$:
\begin{equation}
J_{n,k+1}\geq (-1)^n(P_n...P_1)M_{n-1,k}:(-p+\beta)^*
\label{JgeqM}
\end{equation}
This may be proven in each case by noting the noting the similar 
starting values of $J_{n,2}=(-1)^n(P_n...P_1)p_{n-1,2}p_{n,2}
=p_{n,2}(-1)^n(P_n...P_1)M_{n-1,1}$, and how the former must increase 
before the latter.
We can then use this to get an expression for the numerator $X^{(3)}_{n,k+1}$
from a linear combination of the numerators in the sum.
The second step is add the remaining terms (V) and (VII) to our sum, i.e.
$$
\frac{\chi_{n-1}^3X^{(3)}_{n,k+1}}{R^{(3)}_{n,k+1}}+(V)+ (VII)
=\frac{\chi_{n-1}^4\tilde W_{n,k+1}}{L_{n,k+1}}
$$
and then use
$$\hbox{LCM}_{n,k}=\hbox{LCM}\{R^{(3)}_{n,k+1},J_{n,k+1},(-1)^n(P_n...P_1)M_{n-1,k}\}$$
to determine individual powers in $L_{n,k+1}:(-p+\beta)^*$ by cancelling 
at increases in $\hbox{LCM}_{n,k}:(-p+\beta)^+$;
in a consecutive pair for $p\in {\cal C}^{(n)}_2$, ${\cal C}^{(n)}_3$ or once for 
$p\in {\cal C}^{(n)}_1$. Note that this forms a similar (diagonal) pattern of cancellations
in the major factors of $L_{n,k+1}$ to that seen in $M_{n,k}$ for 
$p\in {\cal C}^{(n)}_1$, ${\cal C}^{(n)}_2$, though in this case the pattern for
$p\in {\cal C}^{(n)}_3$ is similar to that for $p\in {\cal C}^{(n)}_2$.
Also note that as before, all cancellations and subsequent increases will be periodic in
$n$ and $k$, e.g. the cancellations at $L_{4,2},\ L_{4,3}: (-6+\beta)^2[*]$ 
propagates to $L_{10,2},\ L_{10,3};\ L_{10,5},\ L_{10,6}: (-6+\beta)^*[*]$ and so on to 
$L_{16,2},\ L_{16,3};\ L_{16,5},\ L_{16,6};\ L_{16,11},\ L_{16,12}: (-6+\beta)^*[*]$
etc.

Also note that in the case $k=n+1$ the general expansion of $u_{n,n+2}$ above is 
simplified by the fact that
terms (IV)-(VII) drop out, leaving:
$$
u_{n,n+2}=\frac{\alpha^{3n}\beta}{\chi_{n-1}\chi_n^2}.\frac{1}{\chi^4_{n-1}}
\bigg\{
\frac{3\psi_n^2\tilde S_{n-1,n}}{Q_n^2K_{n-1,n}} + \frac{2\psi_n\chi^2_n \tilde W_{n-1,n+1}}{Q_nP_nL_{n-1,n+1}}
+ \frac{\chi^4_n\tilde G_{n-1,n}}{P^2_nJ_{n-1,n}}\bigg\}
$$
This can be further simplified by using the rules for $J_{n,n+1}$, $L_{n,n+2}$ and $K_{n,n+1}$
above to get:
$$
u_{n,n+2}=\frac{\alpha^{3n}\beta^2}{L_{n,n+2}\chi_{n-1}\chi^2_n}.\frac{1}{p^2_{n-1,2}\chi^4_{n-1}}
\bigg\{
3p^2_{n,1}\Omega_{n-1}^2\psi_n^2\tilde S'_{n-1,n}
$$
$$
-2p_{n,1}p_{n,2}\Omega_{n-1}\psi_n\chi^2_n \tilde W'_{n-1,n+1}
+ p_{n,2}^2\chi^4_n\Omega_{n-2}\tilde G'_{n-1,n}
\bigg\}
$$
Hence by using (\ref{Gdnnp1}) we see that can be factorised, leading us to implicitly define 
$\Omega_n$ as:
\begin{equation}
-p_{n-1,2}\chi^2_{n-1}\Omega_n = p_{n,1}\Omega_{n-1}\psi_n + p_{n,2}\Omega_{n-2}\chi^2_n
\label{omegan}
\end{equation}
and $\tilde W'_{n,n+2}$ as:
\begin{equation}
-p_{n-1,2}\chi^2_{n-1} \tilde W'_{n,n+2} = 3p_{n,1}\Omega_{n-1}\psi_n \tilde S'_{n-1,n}
+ p_{n,2}\chi^2_n\tilde G'_{n-1,n},
\label{Wnnp2}
\end{equation}
with which the factorisation of $\tilde W_{n,n+2}=\beta\Omega_n\tilde W'_{n,n+2}$ 
in Postulate \ref{prop1} becomes apparent. We can also get an explicit formula for
$\tilde W_{n,n+2}$:
\begin{equation}
\tilde W'_{n,n+2}=\Omega_{n}\tilde S'_{n-1,n}-2\Omega_{n-1}\tilde S'_{n,n+1}
\label{Wdnnp2}
\end{equation}
We can also combine the formulas (\ref{Wnnp2}) and (\ref{Gnnp1imp}) and use induction 
to show that:
\begin{equation} 
\tilde W'_{n,n+2}+\tilde G'_{n,n+1}=(-1)^np_{n,2}\chi_n^2
\label{WGeqn}
\end{equation}

\subsection{Calculation of $\tau_{n,k}$}
\label{ssKalg}
To simplify the calculation of $\tau_{n,k}$ (and deal with edge cases), it is helpful 
work in terms of {\em modified} numerators:.
$$
\tau_{n,k}=\frac{c_1^{2n+k}\beta \tilde S_{n,k}}{K_{n,k}\chi^{3}_n},\ 
\tilde S_{n,k}=\chi_n^{3-\min\{k,3\}}S_{n,k}
$$
Again using (\ref{aaa}) to cancel the $\chi_{n-2}^3$ terms,
and expanding $x^{2n+k}:T_n$ from (\ref{Tna})
using the postulate, we find:
$$  
\tau_{n,k}=
\frac{\alpha^{2n+k}\beta}{\chi_{n}^3}.\frac{1}{\chi^6_{n-1}}\bigg\{
\bigg(\frac{\psi_n^3\tilde S_{n-1,k-1}}{Q_n^3K_{n-1,k-1}}\bigg)
+
\bigg(\frac{\psi_n^2\tilde W_{n-1,k}\chi_n^{2}
}{Q_n^2P_nL_{n-1,k}}\bigg)
$$
$$+
\bigg(\frac{
\psi_n\tilde G^{(n-1)}_{k-1}\chi_n^4}{Q_nP_n^2J_{n-1,k-1}}\bigg)
+
\bigg(\frac{\tilde S_{n-2,k-2}\chi_n^{6}}
{P_n^3K_{n-2,k-2}}\bigg)
$$
$$
+\bigg(\frac{\tilde W_{n,k+2}
\chi_{n-1}^{5}\chi_n^1}{L_{n,k+2}}\bigg)
+\bigg(\frac{\tilde S_{n-1,k+2}
\chi_{n-1}^{3}\chi_n^3}{K_{n-1,k+2}}\bigg)
- \bigg(\frac{\chi_{n-1}^4\chi_n^2\tilde G^{(n)}_{k+2}}{J_{n,k+2}}\bigg)
$$
$$
-
\bigg(\frac{(-1)^{n-2}\psi_n\varphi^{(n-1)}_k
\chi_{n-1}^{3}\chi_n^{3}}
{(P_n...P_1)Q_nM_{n-1,k}}\bigg)
-
\bigg(\frac{2(-1)^{n-2}\varphi^{(n-1)}_{k+1}
\chi_{n-1}^{4}\chi_n^{4}}
{(P_n...P_1)M_{n-1,k+1}}\bigg)
$$
$$+
\bigg(\frac{2(-1)^{n-3}\psi_n\varphi^{(n-2)}_{k-2}
\chi_{n-1}^{2}\chi_n^{4}
}{Q_nP_n(P_{n}...P_1)M_{n-2,k-2}}\bigg)
+
\bigg(\frac{(-1)^{n-2}\psi_n\varphi^{(n)}_{k}
\chi_{n-1}^{4}\chi_n^{2}
}{Q_n(P_{n}...P_1)M_{n,k}}\bigg)
\bigg\}
$$
We label the components inside the brackets as (I)-(XI) for reference.
As with $u_{n,k+1}$, we can do this sum into two steps; first\footnote{As in the previous case,
this could be decomposed further, i.e. adding (I)+(II)+(III)+(IV) to get to $\chi_{n-1}^2$
before adding remaining terms in sequence to get to $\chi_{n-1}^5$. However, details are omitted
here for the same reason as before.} adding all terms apart from (V)
to get a $\chi_{n-1}^5$ term: 
$$
(I)+(II)+(III)+(IV)+(X)+(VI)+(VIII)+(VII)+(IX)+(XI)
=\frac{\chi_{n-1}^5Y^{(5)}_{n,k}}{V^{(5)}_{n,k}}
$$
where $V^{(5)}_{n,k}$ is determined from the LCM of denominators in this sum,
reducing the power of $(-p+\beta)^*$ in the LCM to match the leading term in 
the denominator of (V) (i.e. $L_{n,k+2}$) for $p\in {\cal C}^{(n)}_1$, ${\cal C}^{(n)}_2$
and ${\cal C}^{(n)}_3$.
We can then use this to get an expression for the numerator $Y^{(5)}_{n,k}$
from a linear combination of the numerators in the sum.
The second step is add the remaining term (V) to our sum, i.e.
$$
\frac{\chi_{n-1}^5Y^{(5)}_{n,k}}{V^{(5)}_{n,k}}+(V)
=\frac{\chi_{n-1}^6\tilde S_{n,k}}{K_{n,k}}
$$
where we use
$$\hbox{LCM}_{n,k}=\hbox{LCM}\{V^{(5)}_{n,k},L_{n,k+2}\}$$
to determine individual powers in $K_{n,k}:(-p+\beta)^*$ by cancelling 
at increases in $\hbox{LCM}_{n,k}:(-p+\beta)^+$;
in a consecutive pair for $p\in {\cal C}^{(n)}_2$, ${\cal C}^{(n)}_3$ or once for 
$p\in {\cal C}^{(n)}_1$.
Note that this forms a similar (diagonal) pattern of cancellations
in the major factors of $K_{n,k}$ to that seen in $L_{n,k+1}$ for 
$p\in {\cal C}^{(n)}_1$, ${\cal C}^{(n)}_2$.
Also note that as before, all cancellations and subsequent increases will be periodic in
$n$ and $k$ (after the second cancel/increase), e.g. the cancellations at 
$K_{4,2},\ K_{4,3}: (-6+\beta)^2[*]$ 
propagates to $K_{10,2},\ K_{10,3};\ K_{10,5},\ K_{10,6}: (-6+\beta)^*[*]$ and so 
on to $K_{16,2},\ K_{16,3};\ K_{16,5},\ K_{16,6};\ K_{16,11},\ K_{16,12}
: (-6+\beta)^*[*]$
etc.

In the case $k=n+1$, all the terms in the general expansion for $\tau_{n,n+1}$ apart from (I)-(IV)
drop out to leave:
$$
\tau_{n,n+1}=\frac{c_1^{3n+1}\beta^2}{\chi^3_n}.\frac{1}{
\chi^6_{n-1}}
\bigg\{
\frac{\psi_n^3\Omega_{n-1}^2\tilde S'_{n-1,n}}{Q_n^3K_{n-1,n}} $$
$$
+ \frac{\psi_n^2\chi^2_n\Omega_{n-1}\tilde W'_{n-1,n+1}}{Q_n^2P_nL_{n-1,n+1}}
+ \frac{\psi_n\chi^4_n\Omega_{n-2}\tilde G'_{n-1,n}}{Q_nP^2_nJ_{n-1,n}}
+ \frac{\chi^6_n\Omega_{n-2}^2\tilde S'_{n-2,n-1}}{P^3_nK_{n-2,n-1}}
\bigg\}
$$
Using the rules above, this may be rewritten as:
$$
\tau_{n,n+1}=\frac{c_1^{3n+1}\beta^2}{K_{n,n+1}\chi^3_n}.\frac{1}{
p_{n-1,2}^3\chi^6_{n-1}}
\bigg\{
p_{n,1}^3\psi_n^3\Omega_{n-1}^2\tilde S'_{n-1,n} 
- p_{n,1}^2p_{n,2}\psi_n^2\chi^2_n\Omega_{n-1}\tilde W'_{n-1,n+1}
$$
$$
+ p_{n,1}p_{n,2}^2\psi_n\chi^4_n\Omega_{n-2}\tilde G'_{n-1,n}
- p_{n,2}^3\chi^6_n\Omega_{n-2}^2\tilde S'_{n-2,n-1}
\bigg\}
$$
Hence using (\ref{Gdnnp1}) and the similar formula for $\tilde W'_{n,n+2}$:
we can factorise the expression, leading us to define $\tilde S'_{n,n+1}$ implicitly by:
\begin{equation}
p_{n-1,2}\chi^2_{n-1}\tilde S'_{n,n+1} = p_{n,1}\psi_n \tilde S'_{n-1,n} - 
p_{n,2}\chi^2_n \tilde S'_{n-2,n-1},
\label{Snnp1}
\end{equation}
for $n\geq 2$, and this gives the corresponding factorisation in 
Postulate~\ref{prop1}.
Also note that since (\ref{phinneqn}) has the same form as (\ref{Snnp1}) and
$S'_{1,2}=\varphi^{(1)}_1=-2$ we can conclude that $\varphi^{(n)}_n=S'_{n,n+1}$ for all $n$.

Finally, considering the general formula for $\tau_{n,k}$ for $k=n$, 
we can see that four terms drop out, i.e.
(V), (VI), (VII) and (VIII). This allows us to write $\tau_{n,n}$ as: 
\begin{equation}
\tau_{n,n}=u_{n,n+2}
+\frac{c_1^{3n}\beta}{K_{n,n}\chi_{n}^3}\bigg\{
\frac{a_3\psi_n^3+a_2\psi_n^2+a_1\psi_n+a_0}{p_{n-1,2}^3(-(n+1)+\beta)\chi^6_{n-1}}
\bigg\}
\label{taunna}
\end{equation}
where
$$
a_3=p_{n,1}^3(-n+b)\tilde S_{n-1,n-1},
$$
$$
a_2=-p_{n,2}p_{n,1}^2\tilde W_{n-1,n}\chi_n^2,
$$
$$
a_1=p_{n,1}p_{n,2}\chi_n^2\big\{
p_{n,2}\chi_n^2\big[\tilde G_{n-1,n-1}+2p_{n-1,2}\varphi^{(n-2)}_{n-2}\chi_{n-1}^2\big]+
p_{n-1,2}^2\varphi^{(n)}_n\chi_{n-1}^4\big\}
$$
$$
a_0=-p_{n,2}^3(-(n-1)+\beta)\tilde S_{n-2,n-2}\chi_n^3
$$
This can then be factored using (\ref{omegan}):
\begin{equation}
\tau_{n,n}=u_{n,n+2}
+\frac{c_1^{3n}\beta}{K_{n,n}\chi_{n}^3}\bigg\{
\frac{\Omega_n(b_2\psi_n^2+b_1\psi_n+b_0)}{p_{n-1,2}^2(-(n+1)+\beta)\chi_{n-1}^4}
\bigg\}
\label{taunnb1}
\end{equation}
where
$$
b_2=p_{n,1}^2(-n+b)\tilde S'_{n-1,n-1},
$$
$$
b_1=p_{n,1}p_{n,2}\chi_n^2\{
[\tilde W_{n-1,n}+(-n+\beta)\Omega_{n-2}\tilde S'_{n-1,n-1}]
-p_{n-1,2}\chi_{n-1}^2\tilde S'_{n-1,n}
\}/\Omega_{n-1}
$$
$$
b_0=p_{n,2}^2(-(n-1)+\beta)\tilde S'_{n-2,n-2}\chi_n^2,
$$
Note that the numerator in $b_1$ is indeed divisible by $\Omega_{n-1}$, leading to various
equalities relating the $\chi_{n-1}$ and $\Omega_{n-1}$ factors. Finally, noting that $u_{n,n+2}$
also has a factor of $\Omega_n$, we factor and combining the terms 
on the right hand side of (\ref{taunnb1}) to get an expression for $\tilde S'_{n,n}$.

\section{Conclusion}
The univariate Pad\'e approximants generated as in Postulate \ref{prop1} will satisfy the 
approximate equation (\ref{Tndefn}) and initial condition $y_n(0)=0$, but not the extra 
condition $y_n(1)=0$. To satisfy this extra condition would require the extension to bivariate
Pad\'e approximants in $x$ and $x^\beta$, similar to what was done for the first order 
Riccati equation in Hegarty, \cite{hegarty}. It would also be interesting to generalise
the methods here to allow for arbitrary polynomial coefficients in the equation and extensions
to the second order Riccati equation from Fair, Luke, \cite{fairluke}, or other 
special cases thereof. This paper will also serve as the foundation for subsequent
work on bivariate Pad\'e approximants for the same equation in \cite{hegarty1}, 
again using the $\tau$-method.

\section*{Acknowledgments}

This work originated during a Post-Doctoral Fellowship in the School of Pharmacy
at the University of Otago under the supervision of Professor Stephen Duffull.

\begin{appendix}
\section{Coefficients}
\label{appA}
\begin{table}
\centering
\begin{tabular}{|c|c|c|}
\hline
$n$ & $p_{n,1}$ & $p_{n,2}$ \\
\hline
$1$ & $1$ & $(-2+\beta)$ \\ 
$2$ & $(-3+\beta)$ & $(-4+\beta)$ \\ 
$3$ & $(-5+\beta)$ & $(-6+\beta)(-2+\beta)$ \\ 
$4$ & $(-7+\beta)$ & $(-8+\beta)$ \\ 
$5$ & $(-9+\beta)(-3+\beta)$ & $(-10+\beta)(-2+\beta)$ \\ 
$6$ & $(-11+\beta)$ & $(-12+\beta)(-4+\beta)$ \\ 
$7$ & $(-13+\beta)$ & $(-14+\beta)(-2+\beta)$ \\ 
$8$ & $(-15+\beta)(-5+\beta)(-3+\beta)$ & $(-16+\beta)$ \\ 
$9$ & $(-17+\beta)$ & $(-18+\beta)(-6+\beta)(-2+\beta)$ \\ 
$10$ & $(-19+\beta)$ & $(-20+\beta)(-4+\beta)$ \\ 
$11$ & $(-21+\beta)(-7+\beta)(-3+\beta)$ & $(-22+\beta)(-2+\beta)$ \\ 
$12$ & $(-23+\beta)$ & $(-24+\beta)(-8+\beta)$ \\ 
\vdots&\vdots&\vdots\\
\hline
\end{tabular} 
\caption{Coefficients $p_{n,1}$ and $p_{n,2}$}
\label{tab1}
\end{table}

\end{appendix}

\end{document}